\documentclass[reqno,11pt]{article}
\usepackage{mathrsfs}
\usepackage{amssymb}
\usepackage[usenames,dvipsnames]{xcolor}
\usepackage{amsthm}
\usepackage{amsmath}
\usepackage{amsfonts}
\usepackage{enumerate}
\usepackage{txfonts}
\usepackage{bm}

\usepackage{graphicx}

\numberwithin{equation}{section}

\newtheorem{theorem}{Theorem}[section]
\newtheorem{lemma}[theorem]{Lemma}
\newtheorem{corollary}[theorem]{Corollary}
\newtheorem{proposition}[theorem]{Proposition}

\theoremstyle{definition}
\newtheorem{definition}[theorem]{Definition}
\newtheorem{assumption}[theorem]{Assumption}
\newtheorem{example}[theorem]{Example}

\theoremstyle{remark}
\newtheorem{remark}[theorem]{Remark}

\begin{document}

\def\be{\begin{eqnarray}}
\def\ee{\end{eqnarray}}
\def\p{\partial}
\def\no{\nonumber}
\def\e{\epsilon}
\def\de{\delta}
\def\De{\Delta}
\def\om{\omega}
\def\Om{\Omega}
\def\f{\frac}
\def\th{\theta}
\def\la{\lambda}
\def\var{\varphi}
\def\b{\bigg}
\def\al{\alpha}
\def\La{\Lambda}
\def\ga{\gamma}
\def\Ga{\Gamma}
\def\ti{\tilde}
\def\Th{\Theta}
\def\si{\sigma}
\def\Si{\Sigma}
\def\wt{\widetilde}
\def\wh{\widehat}
\def\ol{\overline}
\def\ul{\underline}
\def\oo{\infty}
\def\na{\nabla}
\def\q{\quad}
\def\z{\zeta}
\def\co{\coloneqq}
\def\eqq{\eqqcolon}
\def\lab{\label}
\def\ka{\kappa}
\def\bt{\begin{theorem}}
\def\et{\end{theorem}}
\def\bc{\begin{corollary}}
\def\ec{\end{corollary}}
\def\bl{\begin{lemma}}
\def\el{\end{lemma}}
\def\bp{\begin{proposition}}
\def\ep{\end{proposition}}
\def\br{\begin{remark}}
\def\er{\end{remark}}
\def\bd{\begin{definition}}
\def\ed{\end{definition}}
\def\bpf{\begin{proof}}
\def\epf{\end{proof}}
\def\bex{\begin{example}}
\def\eex{\end{example}}
\def\bq{\begin{question}}
\def\eq{\end{question}}
\def\bas{\begin{assumption}}
\def\eas{\end{assumption}}
\def\ber{\begin{exercise}}
\def\eer{\end{exercise}}
\def\mb{\mathbb}
\def\mbR{\mathbb{R}}
\def\mbZ{\mathbb{Z}}
\def\mc{\mathcal}
\def\mcS{\mc{S}}
\def\ms{\mathscr}
\def\lan{\langle}
\def\ran{\rangle}
\def\lb{\llbracket}
\def\rb{\rrbracket}

\title{Existence of axially symmetric weak solutions to steady MHD with non-homogeneous boundary conditions}

\author{Shangkun Weng\thanks{ Pohang Mathematics Institute, Pohang University of Science and Technology. Pohang, Gyungbuk, 790-784, Republic of Korea. Email: skwengmath@gmail.com.}}

\pagestyle{myheadings} \markboth{Existence of axially symmetric weak solutions to steady MHD with non-homogeneous boundary conditions}{Existence of axially symmetric weak solutions}\maketitle

\begin{abstract}
  We establish the existence of axially symmetric weak solutions to steady incompressible magnetohydrodynamics with non-homogeneous boundary conditions. The key issue is the Bernoulli's law for the total head pressure $\Phi=\f 12(|{\bf u}|^2+|{\bf h}|^2)+p$ to a special class of solutions to the inviscid, non-resistive MHD system, where the magnetic field only contains the swirl component. 
\end{abstract}

\begin{center}
\begin{minipage}{5.5in}
\end{minipage}
\end{center}

\section{Introduction and main results}

Let $\Om\subset \mbR^3$ be an axially symmetric domain with $C^2$-smooth boundary $\p\Om= \bigcup_{j=0}^N \Ga_j$ consisting of $N+1$ disjoint components $\Ga_j$; i.e.,
\be\lab{a21}
\Om=\Om_0 \setminus (\cup_{j=1}^N \ol{\Om_j}),\q \ol{\Om_j}\subset \Om_0, j=1,\cdots,N,
\ee
where $\Ga_j=\p\Om_j$. Consider the steady magnetohydrodynamics (MHD) equations in $\Om$:
\be\lab{smhd}\begin{cases}
({\bf u}\cdot\nabla) {\bf u}+\nabla p =({\bf h}\cdot\nabla) {\bf h}+ \Delta {\bf u}+\na\times {\bf f},\q\q\forall x\in \Om,\\
({\bf u}\cdot\nabla){\bf h}-({\bf h}\cdot\nabla){\bf u}=\Delta {\bf h}+ \nabla\times {\bf g},\q\q\forall x\in \Om,\\
\text{div } {\bf u}=\text{div } {\bf h}=0,\q\q\forall x\in \Om,\\
{\bf u} ={\bf a}, \q{\bf h}= {\bf b}\q\q \text{on }\p\Om.
\end{cases}
\ee
For the existence of weak solutions to (\ref{smhd}), the following compatibility conditions are necessary:
\be\lab{comp1}
\sum_{j=0}^N \mc{F}_j\co \sum_{j=0}^N \int_{\Ga_j} {\bf a}\cdot {\bf n} ds=0,\\\lab{comp2}
\sum_{j=0}^N \mc{G}_j\co \sum_{j=0}^N \int_{\Ga_j} {\bf b}\cdot {\bf n} ds=0,
\ee
where ${\bf n}$ is the outward unit vector to the boundary $\p\Om$.

If the magnetic field ${\bf h}$ is absent, then (\ref{smhd}) is reduced to the famous steady Navier-Stokes equations
\be\lab{ns}\begin{cases}
({\bf u}\cdot\nabla) {\bf u}+\nabla p =\Delta {\bf u}+ \na\times {\bf f},\q\q\forall x\in \Om,\\
\text{div } {\bf u}=0,\\
{\bf u} ={\bf a} \q \text{on }\p\Om.
\end{cases}
\ee
Leray \cite{leray33} made fundamental contributions to the existence theory and showed the existence of a weak solution ${\bf u}\in W^{1,2}(\Om)$ to (\ref{ns}) under the stronger assumptions
\be\lab{comp3}
\mc{F}_j= \int_{\Ga_j} {\bf a}\cdot {\bf n} dS=0,\q j=0,1,\cdots, N.
\ee
Leray provided two different methods for the existence results in \cite{leray33}. The first one reduced the nonhomogeneous case to homogeneous case by using the solenoidal extension of boundary value ${\bf a}$ into $\Om$, which was successively completed and clarified in \cite{finn61,hopf,lady69}). The second one is based on a clever contradiction argument, which was used in \cite{amick84,bp94,kp83,rs08}. However, the problem that whether (\ref{ns}), (\ref{comp1}) admit a solution or not is open for long times and usually referred as {\it Leray's problem} in literatures. For sufficiently small fluxes $\mc{F}_j$, one can also obtain the existence of weak solutions \cite{bp94,finn61,fujita59,galdi91,galdi11,kp83,ky09,rs08}. The existence was also known with certain symmetric restrictions on the domain and the boundary data and the forcing term (see \cite{amick84,fujita97,kpr14jmpa,mf02,morimoto95,pr12}). Recently, Korobkov, Pileckas and Russo have made important breakthrough in a series of papers \cite{kpr13arma, kpr15as,kpr15annals,kpr14exterior} on the existence theory without any restrictions on the fluxes. First, in \cite{kpr13arma} they obtained the existence for a plane domain $\Om$ with two connected components of the boundary assuming only the inflow condition on the external component. The new ingredients of analysis in \cite{kpr13arma} are the weak one-sided maximum principle for the total head pressure $\Phi=\f{1}{2}|{\bf u}|^2 +p$ obtained by the Bernoulli's law for weak solutions to the Euler equations and a divergence form representation of $\Phi$. The Bernoulli's law is based on the Morse-Sard theorem developed in \cite{bkk13}. The spatial axially symmetric case was investigated in \cite{kpr15as}, where the existence was established without any restrictions on the fluxes, if all components $\Ga_j$ of $\p\Om$ intersect the axis of symmetry.

In \cite{kpr15annals}, Korobkov, Pileckas and Russo finally established the existence of weak solutions ${\bf u}\in H^1(\Om)$ to the steady Navier-Stokes with boundary values ${\bf a}\in W^{3/2,2}(\p\Om)$ and the force $\na\times {\bf f}\in H^1(\Om)$ in 2-D bounded domain or 3-D axially symmetric domain with $C^2$-smooth boundary, assuming only the total fluxes are zero. By the Morse-Sard theorem proved in \cite{bkk13}, almost all level sets of the stream function $\psi$ are finite unions of $C^1$ curves. Based on the clear understanding of the level sets of $\psi$ and $\Phi$, they can construct appropriate integration domains (bounded by smooth level lines) and estimate the upper bound of the $L^2$ of $\na\Phi$. On the other hand, the length of each of these level lines is bounded from below and the coarea formula implies a lower bound for the $L^1$ norm of $\na \Phi$, from which they can derive a contradiction. In the proof given in \cite{kpr15annals}, {\it the Bernoulli's law} for the Euler equations plays an essential role.


In this paper, we adapt their idea in \cite{kpr15annals} to the steady MHD equations. More precisely, we will establish the existence of axially symmetric weak solutions ${\bf u}({\bf x})= u_r(r,z) {\bf e}_r + u_{\th}(r,z) {\bf e}_{\th} + u_z(r,z){\bf e}_z$ and ${\bf h}({\bf x})= h_{\th}(r,z){\bf e}_{\th}$ to (\ref{smhd}) with nonhomogeneous boundary values in axially symmetric domains with $C^2$ smooth boundary. We introduce some notations. Let $O_{x_1}, O_{x_2}, O_{x_3}$ be coordinate axes in $\mbR^3$ and $\theta =\arctan(x_2/x_1), r=(x_1^2+x_2^2)^{1/2}, z=x_3$ be cylindrical coordinates. Denote by $v_{\th}, v_r, v_z$ the projections of the vector ${\bf v}$ on the axes $\th,r,z$. A function $f$ is said to be {\it axially symmetric} if it does not depend on $\th$. A vector-valued function ${\bf h}= (h_r, h_{\th}, h_z)$ is called {\it axially symmetric} if $h_r, h_{\th}$ and $h_z$ do not depend on $\th$. A vector-valued function ${\bf h}'= (h_r, h_{\th}, h_{z})$ is called {\it axially symmetric with no swirl} if $h_{\th}=0$ while $h_r$ and $h_z$ do not depend on $\th$.

We need to use the following symmetry assumptions.
\begin{description}
  \item[(SO)] $\Om\subset \mbR^3$ is a bounded domain with $C^2$ boundary and $O_{x_3}$ is a symmetry axis of $\Om$.
  \item[(AS)] The assumptions ({\bf SO}) are fulfilled and both the boundary value ${\bf a}\in W^{3/2,2}(\p\Om)$ and $\na\times {\bf f}\in W^{1,2}(\Om)$ are axially symmetric.
  \item[(ASwR)] The assumptions ({\bf SO}) are fulfilled and both the boundary value ${\bf a}\in W^{3/2,2}(\p\Om)$ and $\na\times {\bf f}\in W^{1,2}(\Om)$ are axially symmetric without rotation.
  \item[(ASoS)] The assumptions ({\bf SO}) are fulfilled and both the boundary value ${\bf b}\in W^{3/2,2}(\p\Om)$ and $\na\times {\bf g}\in W^{1,2}(\Om)$ are axially symmetric with only swirl component.
\end{description}
We will use standard notation for Sobolev spaces: $W^{k,q}(\Om), W_0^{k,q}(\Om), W^{\al,q}(\p\Om)$, where $\al\in (0,1), k\in \mb{N}_0, q\in [1,\oo]$. Denote by $H(\Om)$ the subspace of all solenoidal vector fields from $W_0^{1,2}(\Om)$ equipped with the norm $\|{\bf u}\|_{H(\Om)}= \|\na {\bf u}\|_{L^2(\Om)}$. Denote by $L_{AS}^q(\Om)$ ($L_{ASwR}^q(\Om)$) the space of all axially symmetric vector-valued functions (without rotation) in $L^q(\Om)$. Similarly define the spaces $L_{ASoS}^q(\Om)$, $H_{AS}(\Om), H_{ASwR}(\Om), H_{ASoS}^q(\Om)$, $W_{AS}^{1,2}(\Om), W^{1,2}_{ASwR}(\Om), W^{1,2}_{ASoS}(\Om)$, $W_{AS}^{3/2,2}(\p\Om), W^{3/2,2}_{ASwR}(\p\Om), W^{3/2,2}_{ASoS}(\p\Om)$ etc. We denote by $\mc{H}^1$ the one-dimensional Hausdorff measure, i.e., $\mc{H}^1(F)=\lim_{t\to 0+} \mc{H}^1_t(F)$, where
\be\no
\mc{H}_t^1(F) =\inf\b\{\sum_{i=1}^{\oo} \text{diam }F_i: \text{diam }F_i\leq t, F\subset \bigcup_{i=1}^{\oo} F_i\b\}.
\ee

The main result of this paper is stated as follows.
\bt\lab{at41}
{\it Assume that $\Om\subset \mbR^3$ is a bounded axially symmetric domain of type (\ref{a21}) with $C^2$-smooth boundary $\p\Om$. If $(\na\times{\bf f},\na\times{\bf g})\in H_{AS}(\Om)\times H_{ASoS}(\Om)$, $({\bf a},{\bf b})\in W_{AS}^{3/2,2}(\p\Om)\times W_{ASoS}^{3/2,2}(\p\Om)$ and ${\bf a}$ satisfy the compatibility condition (\ref{comp1}), then (\ref{smhd}) admits at least one weak axially symmetric solution $({\bf u}, {\bf h})\in H_{AS}(\Om)\times H_{ASoS}(\Om)$. Moreover, if $\na\times {\bf f}\in H_{ASwR}(\Om)$ and ${\bf a}\in W_{ASwR}^{3/2,2}(\p\Om)$ are axially symmetric with no swirl, then (\ref{smhd}) admits at least one weak axially symmetric solution with $({\bf u}, {\bf h})\in H_{ASwR}(\Om)\times H_{ASoS}(\Om)$.
}
\et
\br\lab{comp}
{\it In the case that ${\bf b}=b_{\th}(r,z){\bf e}_{\th}$, (\ref{comp2}) holds automatically since ${\bf e}_{\th}\cdot {\bf n}\equiv 0$ on $\p\Om$.
}
\er

For the stationary MHD equations (\ref{smhd}), we can define the total head pressure $\Phi=\f{1}{2}(|{\bf u}|^2+|{\bf h}|^2)+ p$. Suppose $({\bf u}, {\bf h}, p)$ are a smooth solution to the inviscid, non-resistive MHD system, then we only have
\be\lab{bernoulli-mhd1}
({\bf u}\cdot\na) \Phi= ({\bf h}\cdot\na) ({\bf u}\cdot {\bf h}).
\ee
So even in the two-dimensional case, the right side is not zero in general. In particular, the level sets of the stream function $\psi$ and $\Phi$ do not coincide with each other, the Bernoulli's law is lost. However, if we further restrict ourself to the axially symmetric MHD case with the special solution form ${\bf u}({\bf x})= u_r(r,z) {\bf e}_r + u_{\th}(r,z) {\bf e}_{\th} + u_z(r,z){\bf e}_z$ and ${\bf h}({\bf x})= h_{\th}(r,z){\bf e}_{\th}$, then $({\bf h}\cdot\na)({\bf u}\cdot{\bf h})=\f{h_{\th}}{r}\p_{\th}({\bf u}\cdot {\bf h})\equiv 0$ and {\it the Bernoulli's law} holds
\be\lab{bernoulli-mhd}
({\bf u}\cdot\nabla) \Phi= 0.
\ee
This has been observed in our previous paper \cite{cw15}, where we have used this to prove some Liouville type theorems for the steady MHD equations. Here we will adapt the methods developed in \cite{kpr15annals} to establish the existence of axially weak weak solutions to (\ref{smhd}).

This paper is organised as follows. We first prepare some preliminaries to reduce the existence problem to some uniform estimates needed in Lemma \ref{pl212} and \ref{pl213}. Then in Section \ref{leray argument}, we first run the Leray's {\it reductio ad absurdum} argument for the steady MHD equations. The Bernoulli's law for the inviscid, nonresistive MHD equations is obtained in Section \ref{inviscid mhd}. Finally, we adapt the methods developed in \cite{kpr15annals} to the steady MHD equation to obtain a contradiction.

\section{Preliminaries}
The following lemmas concern the existence of solenoidal extensions of boundary values.
\bl\lab{pl22}
{\it
\begin{enumerate}[(i)]
  \item If ${\bf a}\in W_{AS}^{3/2,2}(\p\Om)$ and (\ref{comp1}) holds, then there exists an axially symmetric solenoidal extension ${\bf A}\in W^{2,2}(\Om)$ of ${\bf a}$ with the estimate
  \be\lab{p21}
  \|{\bf A}\|_{W_{AS}^{2,2}(\Om)}\leq c \|{\bf a}\|_{W_{AS}^{3/2,2}(\p\Om)}.
  \ee
  Moreover, if conditions ({\bf ASwR}) is prescribed, then ${\bf A}$ can be chosen to have zero swirl component.
  \item If ${\bf b}\in W_{ASoS}^{3/2,2}(\p\Om)$ , then there exists a unique vector field ${\bf H}\in W^{2,2}_{ASoS}(\Om)$ such that
  \be\lab{p222}
  \Delta {\bf H}=0\q \text{in }\Om,\q \text{div }{\bf H}=0\q\text{in }\Om,\q {\bf H}={\bf b}\q \text{on }\p\Om.
  \ee
  We also have the estimate
  \be\lab{p22}
  \|{\bf H}\|_{W_{ASoS}^{2,2}(\Om)}\leq c \|{\bf b}\|_{W_{ASoS}^{3/2,2}(\p\Om)}.
  \ee
\end{enumerate}
}
\el

\bpf
The conclusion (i) has been proved in \cite{kpr15as}. (ii) Let ${\bf b}\in W_{ASoS}^{3/2,2}(\p\Om)$. Then there exists a unique vector field ${\bf F}\in W^{2,2}(\p\Om)$ to the Laplace equation
\be\lab{laplace}
\Delta {\bf F}=0\q\q\text{in }\Om,\q\q\q {\bf F}={\bf b}\q \q \text{on }\p\Om.
\ee
By similar arguments as in Lemma 2.2 in \cite{kpr15as}, we can choose ${\bf F}$ to be axially symmetric. By the standard formulas for $\Delta$ in cylindrical coordinate system, one has for ${\bf F}=(F_r, F_{\th}, F_z)$
\be\lab{la1}
\Delta {\bf F}= (\Delta_2-\f1{r^2}) F_r {\bf e}_r + (\Delta_2-\f1{r^2}) F_{\th} {\bf e}_{\th}+ (\Delta_2 F_z) {\bf e}_z=0,
\ee
where $\Delta_2=(\p_r^2+\f1{r}\p_r+\p_z^2)$. Take ${\bf H}=F_{\theta}{\bf e}_{\theta}$. Then ${\bf H}\in W_{ASoS}^{2,2}(\Om)$ and it follows easily from (\ref{la1}) that
\be\no
\Delta {\bf H}=0.
\ee
Since ${\bf b}\in W_{ASoS}^{3/2,2}(\p\Om)$, we still have ${\bf H}={\bf b}$ on $\p\Om$, therefore ${\bf H}={\bf F}$ by uniqueness. That is, $F_r=F_z\equiv 0$, which implies that
\be\no
\text{div }{\bf H}=\text{div }{\bf F}=\p_r F_r+ \f1{r} F_r +\p_z F_z=0
\ee
\epf
\br\lab{extension}
{\it The statement and proof of (ii) were suggested by one of the referees. The author would like to thank him for the important improvement.
}
\er

Given a function ${\bf F}\in L^q(\Om)$ with $q\geq 6/5$, consider the continuous linear functional $H(\Om)\ni \bm{\eta}\mapsto \int_{\Om}{\bf F}\cdot \bm{\eta} dx$. By the Riesz representation theorem, there exists a unique function ${\bf G}\in H(\Om)$ with
\be\no
\int_{\Om} {\bf F}\cdot \bm{\eta} dx= \int_{\Om}\nabla \bm{\eta} \cdot \nabla {\bf G} dx =\lan {\bf G}, \bm{\eta}\ran_{H(\Om)}\q\q \forall \bm{\eta}\in H(\Om).
\ee

Put ${\bf G}= T_0 {\bf F}$. Evidently, $T_0$ is a continuous linear operator from $L^q(\Om)$ to $H(\Om)$. The following lemmas are easily verified.

\bl\lab{pl29}
{\it The operator $T_0: L^{3/2}(\Om)\to H(\Om)$ has the following symmetry properties:
\be\lab{p28}\begin{array}{ll}
\forall {\bf F}\in L_{AS}^{3/2}(\Om)\q\q T_0 {\bf F}\in H_{AS}(\Om),\\
\forall {\bf F}\in L_{ASwR}^{3/2}(\Om)\q\q T_0 {\bf F}\in H_{ASwR}(\Om),\\
\forall {\bf F}\in L_{ASoS}^{3/2}(\Om)\q\q T_0 {\bf F}\in H_{ASoS}(\Om).
\end{array}
\ee
}
\el

\bl\lab{pl210}
{\it The following inclusions are valid:
\be\lab{p210}\begin{array}{ll}
\forall {\bf u}, {\bf v}\in H_{AS}(\Om)\q\q &({\bf u}\cdot\nabla) {\bf v}\in L_{AS}^{3/2}(\Om), \\
\forall {\bf u}, {\bf v}\in H_{ASwR}(\Om)\q\q &({\bf u}\cdot \nabla) {\bf v}\in L_{ASwR}^{3/2}(\Om),\\
\forall {\bf u}\in H_{AS}(\Om), {\bf v}\in H_{ASoS}(\Om)\q\q &({\bf u}\cdot \nabla) {\bf v}-({\bf v}\cdot \nabla) {\bf u}\in L_{ASoS}^{3/2}(\Om),\\
\forall {\bf u}, {\bf v}\in H_{ASoS}(\Om),\q \q &({\bf u}\cdot\na){\bf v}\in L^{3/2}_{ASwR}(\Om).
\end{array}\ee
}\el

Suppose ${\bf a}\in W^{3/2,2}(\p\Om)$ and also the conditions (\ref{comp1}) and ({\bf AS}) (or ({\bf ASwR})) are fulfilled, then we can find a weak axially symmetric solution ${\bf U}\in W^{2,2}(\Om)$ to the Stokes problem in the sense that ${\bf U}-{\bf A}\in H(\Om)\cap W^{2,2}(\Om)$ and the following formula is satisfied by ${\bf U}$:
\be\no
\int_{\Om} \na {\bf U}\cdot \na \bm{\eta} dx = \int_{\Om}(\na\times {\bf f})\cdot \bm{\eta} dx,\q \forall \bm{\eta}\in H(\Om).
\ee
Moreover,
\be\no
\|{\bf U}\|_{W^{2,2}(\Om)}\leq c(\|{\bf a}\|_{W^{3/2,2}(\p\Om)}+\|\na\times {\bf f}\|_{L^{2}(\Om)}).
\ee
Put ${\bf w}= {\bf u}- {\bf U}$ and ${\bf k}={\bf h}-{\bf H}$. Then the problem (\ref{smhd}) is equivalent to
\be\lab{p212}\begin{cases}
-\Delta {\bf w} + ({\bf U}\cdot\nabla) {\bf w} + ({\bf w}\cdot\nabla) {\bf w} + ({\bf w}\cdot\nabla){\bf U} =-\nabla p- ({\bf U}\cdot\nabla) {\bf U}\\
\q\q\q+({\bf H}\cdot \na){\bf k}+({\bf k}\cdot\na){\bf k}+({\bf k}\cdot\na){\bf H}+({\bf H}\cdot\na){\bf H},\q&\text{in }\Om,\\
-\Delta {\bf k}+({\bf U}\cdot\na){\bf k}+({\bf w}\cdot\na){\bf k}+ ({\bf w}\cdot\na){\bf H}-({\bf k}\cdot\na){\bf U}-({\bf k}\cdot\na){\bf w}-({\bf H}\cdot\na){\bf w}= 0\\
\q\q\q-({\bf U}\cdot\na){\bf H}+ ({\bf H}\cdot\na) {\bf U}+ \na\times {\bf g},\q&\text{in }\Om,\\
\text{div }{\bf w}=\text{div }{\bf k}=0\q&\text{in }\Om,\\
{\bf w}={\bf k}=0\q &\text{on }\p\Om.
\end{cases}\ee

By a {\it weak solution} to the problem (\ref{smhd}) we understand functions $({\bf u},{\bf h})$ such that ${\bf w}= {\bf u}- {\bf U}\in H(\Om)$, ${\bf k}={\bf h}-{\bf H}\in H(\Om)$ and for any $\bm{\eta}\in H(\Om), \bm{\z}\in W_0^{1,2}(\Om)$
\be\lab{p213}\begin{array}{ll}
\lan {\bf w}, \bm{\eta}\ran_{H(\Om)} &= -\int_{\Om} ({\bf U}\cdot \nabla) {\bf U}\cdot \bm{\eta} dx -\int_{\Om} ({\bf U}\cdot\na){\bf w}\cdot \bm{\eta} dx -\int_{\Om} ({\bf w}\cdot\nabla){\bf w}\cdot\bm{\eta} dx\\
&\q- \int_{\Om}({\bf w}\cdot\nabla){\bf U}\cdot \bm{\eta} dx+\int_{\Om}({\bf H}\cdot\na){\bf k}\cdot\bm{\eta} dx + \int_{\Om}({\bf k}\cdot\na){\bf k}\cdot \bm{\eta} dx\\&\q + \int_{\Om}({\bf k}\cdot\na){\bf H}\cdot \bm{\eta} dx+ \int_{\Om}({\bf H}\cdot\na){\bf H}\cdot \bm{\eta} dx, \\
\lan {\bf k},\bm{\zeta}\ran_{H(\Om)}&=-\int_{\Om}({\bf U}\cdot\na){\bf H}\cdot\bm{\z} dx-\int_{\Om}({\bf U}\cdot\na){\bf k}\cdot\bm{\z} dx- \int_{\Om}({\bf w}\cdot\na){\bf k}\cdot \bm{\z} dx\\&\q-\int_{\Om}({\bf w}\cdot\na){\bf H}\cdot\bm{\z} dx+\int_{\Om}({\bf k}\cdot\na){\bf U}\cdot\bm{\z} dx+\int_{\Om}({\bf k}\cdot\na){\bf w}\cdot \bm{\zeta}\\&\q+ \int_{\Om}({\bf H}\cdot\na){\bf w}\cdot\bm{\z}dx +\int_{\Om}({\bf H}\cdot\na){\bf U}\cdot\bm{\z} dx+ \int_{\Om}(\na\times {\bf g})\cdot \bm{\z} dx.
\end{array}\ee

By the Riesz representation theorem, for any $\left(\begin{array}{ll}{\bf w}\\{\bf k}\end{array}\right)\in H(\Om)\times H(\Om)$ there exists a unique element $\mb{T}\left(\begin{array}{ll}{\bf w}\\{\bf k}\end{array}\right)=\b(T_1\left(\begin{array}{ll}{\bf w}\\{\bf k}\end{array}\right), T_2\left(\begin{array}{ll}{\bf w}\\{\bf k}\end{array}\right)\b)^T\in H(\Om)\times H(\Om)$ such that the right-hand sides of (\ref{p213}) are equivalent to $\b\lan T_1\left(\begin{array}{ll}{\bf w}\\{\bf k}\end{array}\right), \bm{\eta}\b\ran_{H(\Om)}$ and $\b\lan T_2\left(\begin{array}{ll}{\bf w}\\{\bf k}\end{array}\right), \bm{\z}\b\ran_{H(\Om)}$ for all $\bm{\eta}\in H(\Om),\bm{\z}\in W_0^{1,2}(\Om)$, respectively. Obviously, $\mb{T}$ is a nonlinear operator from $H(\Om)\times H(\Om)$ to $H(\Om)\times H(\Om)$.

\bl\lab{pl211}
{\it The operator $\mb{T}: H(\Om)\times H(\Om)\to H(\Om)\times H(\Om)$ is compact. Moreover, $\mb{T}$ has the following symmetry properties:
\be\lab{p214}\begin{array}{ll}
\forall \left(\begin{array}{ll}{\bf w}\\{\bf k}\end{array}\right)\in H_{AS}(\Om)\times H_{ASoS}(\Om),\q T_1 \left(\begin{array}{ll}{\bf w}\\{\bf k}\end{array}\right)\in H_{AS}(\Om),\\
\forall \left(\begin{array}{ll}{\bf w}\\{\bf k}\end{array}\right)\in H_{ASwR}(\Om)\times H_{ASoS}(\Om),\q T_1 \left(\begin{array}{ll}{\bf w}\\{\bf k}\end{array}\right)\in H_{ASwR}(\Om),\\
\forall \left(\begin{array}{ll}{\bf w}\\{\bf k}\end{array}\right)\in H_{AS}(\Om)\times H_{ASoS}(\Om),\q T_2 \left(\begin{array}{ll}{\bf w}\\{\bf k}\end{array}\right)\in H_{ASoS}(\Om).
\end{array}\ee
}\el
\bpf
The compactness can be proved in a standard way as shown in \cite{lady69} and (\ref{p214}) follows from Lemma \ref{pl29} and Lemma \ref{pl210}.
\epf
Hence (\ref{p213}) is equivalent to the operator equation
\be\lab{p216}
\left(\begin{array}{ll}{\bf w}\\{\bf k}\end{array}\right)=\mb{T}\left(\begin{array}{ll}{\bf w}\\{\bf k}\end{array}\right)
\ee
in the space $H(\Om)\times H(\Om)$. Thus, we can apply the Leray-Schauder fixed point theorem to the compact operators $\mb{T}|_{H_{AS}(\Om)\times H_{ASoS}(\Om)}$ and $\mb{T}|_{H_{ASwR}(\Om)\times H_{ASoS}(\Om)}$. Then the following statements hold.

\bl\lab{pl212}
{\it Let conditions ({\bf AS})-({\bf ASoS}), (\ref{comp1})-(\ref{comp2}) be fulfilled. Suppose all possible solutions $\left(\begin{array}{ll}{\bf w}\\{\bf k}\end{array}\right)$ to the equation $\left(\begin{array}{ll}{\bf w}\\{\bf k}\end{array}\right)= \la \mb{T}\left(\begin{array}{ll}{\bf w}\\{\bf k}\end{array}\right)$ with $\la\in [0,1]$are uniformly bounded in $H(\Om)\times H(\Om)$. Then problem (\ref{smhd}) admits at least one weak axially symmetric solution $({\bf u}, {\bf h})\in H_{AS}(\Om)\times H_{ASoS}(\Om)$.
}
\el

\bl\lab{pl213}
{\it Let conditions ({\bf ASwR})-({\bf ASoS}), (\ref{comp1})-(\ref{comp2}) be fulfilled. Suppose all possible solutions $\left(\begin{array}{ll}{\bf w}\\{\bf k}\end{array}\right)$ to the equation $ \left(\begin{array}{ll}{\bf w}\\{\bf k}\end{array}\right)= \la \mb{T}\left(\begin{array}{ll}{\bf w}\\{\bf k}\end{array}\right)$ with $\la\in [0,1]$ are uniformly bounded in $H(\Om)\times H(\Om)$. Then problem (\ref{smhd}) admits at least one weak axially symmetric solution $({\bf u}, {\bf h})\in H_{ASwR}(\Om)\times H_{ASoS}(\Om)$.
}
\el

\section{Proof of Theorem \ref{at41}}\lab{proof}

\subsection{The {\it reductio ad absurdum} argument by Leray}\lab{leray argument}

We apply the {\it reductio ad absurdum} argument of Leray \cite{leray33} to the stationary MHD equations. To prove the existence of a weak solution to the MHD system (\ref{smhd}), by Lemma \ref{pl212} and \ref{pl213} it is sufficient to show that the weak solutions $({\bf w}, {\bf k})$  satisfying for any $(\bm{\eta}, \bm{\z})\in H(\Om)\times W_0^{1,2}(\Om)$
\be\lab{o44}\begin{array}{ll}
\lan {\bf w}, \bm{\eta}\ran_{H(\Om)} &= -\la\int_{\Om} ({\bf U}\cdot \nabla) {\bf U}\cdot \bm{\eta} dx -\la\int_{\Om} ({\bf U}\cdot\na){\bf w}\cdot \bm{\eta} dx -\la\int_{\Om} ({\bf w}\cdot\nabla){\bf w}\cdot\bm{\eta} dx\\
&\q- \la\int_{\Om}({\bf w}\cdot\nabla){\bf U}\cdot \bm{\eta} dx+\la\int_{\Om}({\bf H}\cdot\na){\bf k}\cdot\bm{\eta} dx +\la \int_{\Om}({\bf k}\cdot\na){\bf k}\cdot \bm{\eta} dx \\
&\q+\la \int_{\Om}({\bf k}\cdot\na){\bf H}\cdot \bm{\eta} dx+ \la\int_{\Om}({\bf H}\cdot\na){\bf H}\cdot \bm{\eta} dx,\\
\lan {\bf k},\bm{\zeta}\ran_{H(\Om)}&=-\la\int_{\Om}({\bf U}\cdot\na){\bf H}\cdot\bm{\z} dx-\la\int_{\Om}({\bf U}\cdot\na){\bf k}\cdot\bm{\z} dx-\la\int_{\Om}({\bf w}\cdot\na){\bf k}\cdot \bm{\z} dx\\&\q-\la\int_{\Om}({\bf w}\cdot\na){\bf H}\cdot\bm{\z} dx+\la\int_{\Om}({\bf k}\cdot\na){\bf U}\cdot\bm{\z} dx+\la\int_{\Om}({\bf k}\cdot\na){\bf w}\cdot \bm{\zeta}\\&\q+\la \int_{\Om}({\bf H}\cdot\na){\bf w}\cdot\bm{\z}dx \la+\int_{\Om}({\bf H}\cdot\na){\bf U}\cdot\bm{\z} dx-\la \int_{\Om}\na {\bf H}\cdot \na\bm{\z} dx + \la\int_{\Om}(\na\times {\bf g})\cdot \bm{\z} dx,
\end{array}\ee
are uniformly bounded in $H(\Om)\times H(\Om)$ with respect to $\la\in [0,1]$. Assume that this is false. Then there exist sequences $\{\la_n\}_{n\in \mb{N}}\subset [0,1]$ and $\{\wh{{\bf w}}_n, \wh{{\bf k}}_n\}_{n\in\mb{N}} \in H(\Om)\times H(\Om)$ such that for any $(\bm{\eta}, \bm{\z})\in H(\Om)\times W_0^{1,2}(\Om)$
\be\lab{o45}\begin{array}{ll}
\int_{\Om} \nabla \wh{{\bf w}}_n\cdot \na\bm{\eta} dx -\la_n\int_{\Om} ((\wh{{\bf w}}_n+ {\bf U})\cdot\nabla) \bm{\eta} \cdot \wh{{\bf w}}_n dx -\la_n\int_{\Om} (\wh{{\bf w}}_n\cdot\nabla) \bm{\eta}\cdot {\bf U} dx\\\q\q+\la_n\int_{\Om}((\wh{{\bf k}}_n+ {\bf H})\cdot\nabla) \bm{\eta} \cdot \wh{{\bf k}}_n dx+\la_n\int_{\Om}(\wh{{\bf k}}_n\cdot\na)\bm{\eta}\cdot {\bf H}dx\\
=\la_n\int_{\Om} ({\bf U}\cdot\nabla) \bm{\eta}\cdot {\bf U} dx-\la_n\int_{\Om} ({\bf H}\cdot\nabla) \bm{\eta}\cdot {\bf H} dx,
\end{array}\\\lab{o45'}\begin{array}{ll}
\int_{\Om}\na \wh{{\bf k}}_n\cdot \na\bm{\z} dx-\la_n\int_{\Om}((\wh{{\bf w}}_n+ {\bf U})\cdot\na)\bm{\z}\cdot \wh{{\bf k}}_n dx-\la_n\int_{\Om}(\wh{{\bf w}}_n\cdot\na)\bm{\z}\cdot {\bf H} dx \\\q\q+\la_n\int_{\Om}((\wh{{\bf k}}_n+{\bf H})\cdot\na)\bm{\z}\cdot \wh{{\bf w}}_n dx+ \la_n \int_{\Om}(\wh{{\bf k}}_n\cdot\na)\bm{\z}\cdot {\bf U} dx\\
=\la_n\int_{\Om}(({\bf U}\cdot\na)\bm{\z})\cdot {\bf H}dx -\la_n\int_{\Om}({\bf H}\cdot\na)\bm{\z}\cdot{\bf U} dx-\la_n\int_{\Om}\na {\bf H}\cdot \na\bm{\z} dx+ \la_n \int_{\Om}(\na\times {\bf g})\cdot \bm{\z} dx
\end{array}\ee
and
\be\lab{o46}
\lim_{n\to\oo} \la_n =\la_0\in [0,1],\q \lim_{n\to\oo} J_n^2 =\lim_{n\to \oo} (\|\wh{{\bf w}}_n\|_{H(\Om)}^2+\|\wh{{\bf k}}_n\|_{H(\Om)}^2)=\oo.
\ee

Denote ${\bf w}_n= J_n^{-1}\wh{{\bf w}}_n, {\bf k}_n= J_n^{-1}\wh{{\bf k}}_n$. Since $\|{\bf w}_n\|_{H(\Om)}^2+\|{\bf k}_n\|_{H(\Om)}^2=1$, there exists a subsequence $\{{\bf w}_{n_l}, {\bf k}_{n_l}\}$ converging weakly in $H(\Om)$ to vector fields ${\bf w},{\bf k}\in H(\Om)$. Because of the compact embedding
\be\no
H(\Om)\mapsto L^r(\Om)\q\q \forall r\in [1,6),
\ee
the subsequence $\{{\bf w}_{n_l}, {\bf k}_{n_l}\}$ converges strongly in $L^r(\Om)$. Replacing $\z$ in (\ref{o45'}) by $J_n^{-2}\z$, and letting $n\to \oo$, we obtain
\be\lab{o48'}
\la_0 \int_{\Om} [({\bf w}\cdot\na){\bf k}-({\bf k}\cdot\na){\bf w}]\cdot \z dx=0.
\ee

Taking $\bm{\eta}= J_n^{-2} \wh{{\bf w}}_n, \bm{\z}= J_n^{-2} \wh{{\bf k}}_n$ in (\ref{o45})-(\ref{o45'}) and adding the above two identities, we get
\be\lab{o47}\begin{array}{ll}
\int_{\Om} |\na {\bf w}_n|^2+|\na{\bf k}_n|^2 dx=\la_n \int_{\Om}[({\bf w}_n\cdot\nabla) {\bf w}_n-({\bf k}_n\cdot){\bf k}_n]\cdot {\bf U} dx-\la_n\int_{\Om}[({\bf w}_n\cdot\na){\bf k}_n-({\bf k}_n\cdot\na){\bf w}_n]\cdot {\bf H} dx\\
\q+J_n^{-1}\la_n \int_{\Om} [({\bf U}\cdot\na){\bf w}_n \cdot {\bf U}-({\bf H}\cdot\na){\bf w}_n\cdot{\bf H}+({\bf U}\cdot\na){\bf k}_n\cdot{\bf H}-({\bf H}\cdot\na){\bf k}_n\cdot {\bf U}] dx\\
\q- J_n^{-1}\la_n\int_{\Om}[(\na\times {\bf g})\cdot {\bf k}_n+\na {\bf H}\cdot\na {\bf k}_n ]dx
\end{array}\ee
 Therefore, passing to a limit as $n_l\to \oo$ in equality (\ref{o47}) and using (\ref{o48'}) we obtain
\be\lab{o48}
1=\la_0\int_{\Om}[({\bf w}\cdot\na){\bf w}-({\bf k}\cdot\na){\bf k}]\cdot {\bf U} dx.
\ee
This implies $\la_0\in (0,1]$. Let us return to the integral identity (\ref{o45}). Consider the functional
\be\no
R_n(\bm{\eta})&=& \int_{\Om} \nabla \wh{{\bf w}}_n\cdot \na\bm{\eta} dx -\la_n\int_{\Om} ((\wh{{\bf w}}_n+ {\bf U})\cdot\nabla) \bm{\eta} \cdot \wh{{\bf w}}_n dx -\la_n\int_{\Om} (\wh{{\bf w}}_n\cdot\nabla) \bm{\eta}\cdot {\bf U} dx\\\no&\q&+\la_n\int_{\Om}((\wh{{\bf k}}_n+ {\bf H})\cdot\nabla) \bm{\eta} \cdot \wh{{\bf k}}_n dx+\la_n\int_{\Om}(\wh{{\bf k}}_n\cdot\na)\bm{\eta}\cdot {\bf H}dx
-\la_n\int_{\Om} ({\bf U}\cdot\nabla) \bm{\eta}\cdot {\bf U} dx\\\no&\q&+\la_n\int_{\Om} ({\bf H}\cdot\nabla) \bm{\eta}\cdot {\bf H} dx\q\q \forall \bm{\eta}\in W_0^{1,2}(\Om).
\ee
Obviously, $R_k(\bm{\eta})$ is a linear functional and
\be\no
|R_n(\bm{\eta})|\leq c(\|(\wh{{\bf w}}_n,\wh{{\bf k}}_n)\|_{H(\Om)}+\|(\wh{{\bf w}}_n,\wh{{\bf k}}_n)\|_{H(\Om)}^2+\|({\bf a},{\bf b})\|_{W^{3/2,2}(\p\Om)}^2+\|{\bf f}\|_{W^{1,2}_0(\Om)}^2)\|\bm{\eta}\|_{H(\Om)}
\ee
with constant $c$ independent of $n$. It follows from (\ref{o45}) that
\be\no
R_n(\bm{\eta})=0\q\q \forall \bm{\eta}\in H(\Om).
\ee
Therefore, there exists an axially symmetric function $\hat{p}_n\in\hat{L}^2(\Om)=\{q\in L^2(\Om): \int_{\Om} q(x) dx=0\}$ such that
\be\no
R_n(\bm{\eta}) =\int_{\Om} \hat{p}_n \text{div }\bm{\eta} dx\q\q \forall \bm{\eta}\in W_0^{1,2}(\Om)
\ee
and
\be\lab{o49}
\|\hat{p}_n\|_{L^2(\Om)}\leq c(\|(\wh{{\bf w}}_n,\wh{{\bf k}}_n)\|_{H(\Om)}+\|(\wh{{\bf w}}_n,\wh{{\bf k}}_n)\|_{H(\Om)}^2+\|({\bf a},{\bf b})\|_{W^{3/2,2}(\p\Om)}^2+\|{\bf f}\|_{W^{1,2}_0(\Om)}^2).
\ee
The pair $(\wh{{\bf w}}_n, \wh{{\bf k}}_n, \hat{p}_n)$ satisfies the integral identity
\be\no
&&\int_{\Om}\nabla \wh{{\bf w}}_n\cdot \na\bm{\eta} dx-\la_n\int_{\Om} ((\wh{{\bf w}}_n+ {\bf U})\cdot\nabla) \bm{\eta} \cdot \wh{{\bf w}}_n dx -\la_n\int_{\Om} (\wh{{\bf w}}_n\cdot\nabla) \bm{\eta}\cdot {\bf U} dx\\\lab{o410}
&&+\la_n\int_{\Om} ((\wh{{\bf k}}_n+ {\bf H})\cdot\nabla) \bm{\eta} \cdot \wh{{\bf k}}_n dx + \la_n\int_{\Om} (\wh{{\bf k}}_n\cdot\nabla) \bm{\eta}\cdot {\bf H} dx+\la_n \int_{\Om} ({\bf H}\cdot\na)\bm{\eta}\cdot {\bf H} dx\\\no
&&\q-\la_n\int_{\Om} ({\bf U}\cdot\nabla) \bm{\eta}\cdot {\bf U} dx=\int_{\Om} \hat{p}_n \text{div }\bm{\eta} dx,\q \forall \bm{\eta}\in W_0^{1,2}(\Om).
\ee
Let $\wh{{\bf u}}_n=\wh{{\bf w}}_n +{\bf U}, \wh{{\bf h}}_n= \wh{{\bf k}}_n+ {\bf H}$. Then identity (\ref{o410}) reduces to
\be\no
&&\int_{\Om}\na \wh{{\bf u}}_n\cdot\nabla \bm{\eta} dx -\int_{\Om} \hat{p}_n\text{div }\bm{\eta} dx= -\la_n \int_{\Om} (\wh{{\bf u}}_n\cdot\na)\wh{{\bf u}}_n \cdot\bm{\eta} dx\\\no
&&\q\q\q\q+\la_n \int_{\Om} (\wh{{\bf h}}_n\cdot\na)\wh{{\bf h}}_n \cdot\bm{\eta} dx+ \la_n \int_{\Om}(\nabla\times {\bf f})\cdot \bm{\eta} dx,\q \forall \bm{\eta}\in W_0^{1,2}(\Om).
\ee

Thus $(\wh{{\bf u}}_n, \wh{{\bf h}}_n, \hat{p}_n)$ might be considered as a weak solution to the Stokes problem
\be\no\begin{cases}
-\De \wh{{\bf u}}_n +\na \hat{p}_n=-\la_n (\wh{{\bf u}}_n\cdot\na) \wh{{\bf u}}_n+ \la_n (\wh{{\bf h}}_n\cdot\na) \wh{{\bf h}}_n+\la_n \na\times {\bf f}\co {\bf F}_n\q &\text{in }\Om,\\
-\De \wh{{\bf h}}_n = -\la_n (\wh{{\bf u}}_n\cdot\na) \wh{{\bf h}}_n+ \la_n (\wh{{\bf h}}_n\cdot\na) \wh{{\bf u}}_n+\na\times {\bf g}\co {\bf H}_n\q &\text{in }\Om,\\
\text{div }\wh{{\bf u}}_n=\text{div }\wh{{\bf h}}_n=0\q &\text{in }\Om,\\
\wh{{\bf u}}_n= {\bf a}, \wh{{\bf h}}_n= {\bf b}\q &\text{on }\p\Om.
\end{cases}\ee
Obviously, ${\bf F}_n, {\bf H}_n\in L^{3/2}(\Om)$ and
\be\no
\|{\bf F}_n\|_{L^{3/2}(\Om)} &\leq& c\|(\wh{{\bf u}}_n\cdot\na)\wh{{\bf u}}_n\|_{L^{3/2}(\Om)}+ c\|(\wh{{\bf h}}_n\cdot\na)\wh{{\bf h}}_n\|_{L^{3/2}(\Om)}+\|\na\times {\bf f}\|_{L^{3/2}(\Om)}\\\no
&\leq& c\|\wh{{\bf u}}_n\|_{L^{6}(\Om)}\|\na \wh{{\bf u}}_n\|_{L^2(\Om)}+ c\|\wh{{\bf h}}_n\|_{L^{6}(\Om)}\|\na \wh{{\bf h}}_n\|_{L^2(\Om)}+\|{\bf f}\|_{W^{1,2}_0(\Om)}\\\no
&\leq& c(\|\wh{{\bf w}}_n\|_{H(\Om)}^2+\|\wh{{\bf k}}_n\|_{H(\Om)}^2+ \|{\bf a}\|_{W^{1/2,2}(\p\Om)}^2+ \|{\bf b}\|_{W^{1/2,2}(\p\Om)}^2)+\|{\bf f}\|_{W^{1,2}_0(\Om)},\\\no
\|{\bf H}_n\|_{L^{3/2}(\Om)} &\leq& c\|(\wh{{\bf u}}_n\cdot\na)\wh{{\bf h}}_n\|_{L^{3/2}(\Om)}+ c\|(\wh{{\bf h}}_n\cdot\na)\wh{{\bf u}}_n\|_{L^{3/2}(\Om)}+\|\na\times {\bf g}\|_{L^{3/2}(\Om)}\\\no
&\leq& c(\|\wh{{\bf w}}_n\|_{H(\Om)}^2+\|\wh{{\bf k}}_n\|_{H(\Om)}^2+ \|{\bf a}\|_{W^{1/2,2}(\p\Om)}^2 +\|{\bf b}\|_{W^{1/2,2}(\p\Om)}^2)+\|{\bf g}\|_{W^{1,2}_0(\Om)},
\ee
where $c$ is independent of $n$. By the well-known regularity results for the Stokes system (see Theorem IV.6.1 in \cite{galdi11}), we have $\wh{{\bf u}}_n, \wh{{\bf h}}_n\in W^{2,3/2}(\Om), \hat{p}_n\in W^{1,3/2}(\Om)$, and also the estimate
\be\lab{a398}\begin{array}{ll}
&\|\wh{{\bf u}}_n\|_{W^{2,3/2}(\Om)}+ \|\hat{p}_n\|_{W^{1,3/2}(\Om)} \leq c(\|{\bf F}_n\|_{L^{3/2}(\Om)}+\|{\bf a}\|_{W^{3/2,2}(\p\Om)})\\
\leq& c(\|\wh{{\bf w}}_n\|_{H(\Om)}^2+\|\wh{{\bf k}}_n\|_{H(\Om)}^2+ \|({\bf a},{\bf b})\|_{W^{3/2,2}(\p\Om)}^2+ \|({\bf a}, {\bf b})\|_{W^{3/2,2}(\p\Om)}+\|{\bf f}\|_{W^{1,2}_0(\Om)}),\end{array}\\\lab{a399}\begin{array}{ll}
&\|\wh{{\bf h}}_n\|_{W^{2,3/2}(\Om)} \leq c(\|{\bf H}_n\|_{L^{3/2}(\Om)}+\|{\bf b}\|_{W^{3/2,2}(\p\Om)})\\
\leq& c(\|\wh{{\bf w}}_n\|_{H(\Om)}^2+\|\wh{{\bf k}}_n\|_{H(\Om)}^2+ \|({\bf a},{\bf b})\|_{W^{3/2,2}(\p\Om)}^2+ \|({\bf a}, {\bf b})\|_{W^{3/2,2}(\p\Om)}+\|{\bf g}\|_{W^{1,2}_0(\Om)}).
\end{array}
\ee

Denote ${\bf u}_n=J_n^{-1}\wh{{\bf u}}_n, {\bf h}_n= J_n^{-1}\wh{{\bf h}}_n$ and $p_n= \la_n^{-1} J_n^{-2} \hat{p}_n$. Then
\be\lab{a310}\begin{array}{ll}
-\nu_n \De {\bf u}_n +({\bf u}_n\cdot\na){\bf u}_n +\na p_n= ({\bf h}_n\cdot\na) {\bf h}_n +\na\times {\bf f}_n,\q&\textit{in }\Om,\\
-\nu_n \De {\bf h}_n +({\bf u}_n\cdot\na){\bf h}_n -({\bf h}_n\cdot\na) {\bf u}_n=\na\times {\bf g}_n,\q&\textit{in }\Om,\\
\text{div }{\bf u}_n= \text{div }{\bf h}_n=0,\q&\textit{in }\Om,\\
{\bf u}_n= {\bf a}_n,\q {\bf h}_n= {\bf b}_n,\q &\textit{on }\p\Om,
\end{array}\ee
where $\nu_n= \la_n^{-1}J_n^{-1}, {\bf f}_n= J_n^{-2} {\bf f}, {\bf g}_n= J_n^{-2} {\bf g}$ and ${\bf a}_n= J_n^{-1}{\bf a}, {\bf b}_n= J_n^{-1}{\bf b}$.

It follows from (\ref{a398}) that
\be\no
\|p_n\|_{W^{1,3/2}(\Om)}\leq \textit{const}.
\ee
Hence, from the sequence $\{p_{n_l}\}$ we can extract a subsequence, still denoted by $\{p_{n_l}\}$, which converges weakly in $W^{1,3/2}(\Om)$ to some function $p\in W^{1,3/2}(\Om)$. Let $\bm{\var}\in C_0^{\oo}(\Om)$. Taking $\bm{\eta} =J_n^{-2} \bm{\var}$ in (\ref{o410}) and letting $n_l\to \oo$, we get
\be\no
-\la_0 \int_{\Om} ({\bf w}\cdot\na) \bm{\var}\cdot {\bf w} dx+\la_0 \int_{\Om} ({\bf k}\cdot \na)\bm{\var}\cdot {\bf k} dx =\la_0\int_{\Om} p \text{div }\bm{\var} dx\q\q\forall \bm{\var}\in C_0^{\oo}(\Om).
\ee
Integrating by parts in the last equality, we derive
\be\lab{o412}
\la_0 \int_{\Om} [({\bf w}\cdot\na){\bf w}-({\bf k}\cdot\na){\bf k}]\cdot \bm{\var} dx = -\la_0\int_{\Om}\na p\cdot \bm{\var} dx\q\forall \bm{\var}\in C_0^{\oo}(\Om).
\ee
Hence, the pair $({\bf w}, {\bf k}, p)$ satisfies, for almost all $x\in\Om$, the inviscid, nonresistive MHD equations
\be\lab{o413}\begin{cases}
({\bf w}\cdot\na) {\bf w}+ \na p=({\bf k}\cdot\na){\bf k},\q &\textit{in }\Om,\\
({\bf w}\cdot\na) {\bf k}- ({\bf k}\cdot\na){\bf w} =0,\q &\textit{in }\Om,\\
\text{div }{\bf w}= \text{div }{\bf k}=0,\q &\textit{in }\Om,\\
{\bf w}={\bf k}=0,\q \q &\textit{on }\p\Om.
\end{cases}\ee

We summarize the above results as follows.
\bl\lab{al41}
{\it Assume that $\Om \subset \mbR^3$ is a bounded axially symmetric domain of type (\ref{a21}) with $C^2$-smooth boundary $\p\Om$, $(\na\times {\bf f},\na\times {\bf g})\in W_{AS}^{1,2}(\Om)\times W_{ASoS}^{1,2}(\Om)$, $({\bf a},{\bf b})\in W_{AS}^{3/2,2}(\p\Om)\times W_{ASoS}^{3/2,2}(\p\Om)$ are axially symmetric, and ${\bf a}$ and ${\bf b}$ satisfy conditions (\ref{comp1})-(\ref{comp2}). If the assertion of Theorem \ref{at41} is false, then there exist ${\bf w}, {\bf k}, p$ with the following properties:
\begin{description}
  \item[(IMHD-AX)] The axially symmetric functions $({\bf w}, {\bf k})\in H_{AS}(\Om)\times H_{ASoS}(\Om), p\in W_{AS}^{1,3/2}(\Om)$ satisfy the invisicd nonresistive MHD system (\ref{o413}) and (\ref{o48}).
  \item[(MHD-AX)] There exist a sequence of axially symmetric functions ${\bf u}_n\in W_{AS}^{1,2}(\Om), {\bf h}_n\in W_{ASoS}^{1,2}(\Om)$, $p_n\in W_{AS}^{1,3/2}(\Om)$ and numbers $\nu_n\to 0+, \la_n\to \la_0\in (0,1]$ such that the norms $\|{\bf u}_n\|_{W^{1,2}(\Om)}+\|{\bf h}_n\|_{W^{1,2}(\Om)}$ and $\|p_n\|_{W^{1,3/2}(\Om)}$ are uniformly bounded, the pair $({\bf u}_n, {\bf h}_n, p_n)$ satisfies (\ref{a310}), and
      \be\lab{a41}\begin{array}{ll}
      \|\nabla {\bf u}_n\|_{L^2(\Om)}+ \|\na {\bf h}_n\|_{L^2(\Om)}\to 1,\\
      {\bf u}_n \rightharpoonup {\bf w},\q {\bf h}_n \rightharpoonup {\bf k}\q \text{in } W^{1,2}(\Om),\q p_n \rightharpoonup p\q \text{in } W^{1,3/2}(\Om).
      \end{array}\ee
      Moreover, $({\bf u}_n, {\bf h}_n)\in W^{3,2}_{loc}(\Om)$ and $p_n\in W^{2,2}_{loc}(\Om)$.
\end{description}
}
\el

Assume that
\be\no
\Ga_j \cap O_{x_3} &\neq & \emptyset, \q j=0,\cdots, M',\\\no
\Ga_j\cap O_{x_3}  &=& \emptyset,\q  j=M'+1,\cdots, N.
\ee

Let $P_+= \{(0,x_2,x_3): x_2>0, x_3\in \mbR\}$, $\mc{D}= \Om\cap P_+$. Obviously, on $P_+$ the coordinates $x_2, x_3$ coincide with the coordinates $r,z$.
For a set $A\subset \mbR^3$, put $\text{\u{A}} \co A\cap P_+$, and for $B\subset P_+$, denote by $\ti{B}$ the set in $\mbR^3$ obtained by rotation of $B$ around the $O_z$-axis. Then
\begin{description}
  \item[($S_1$)] $\mc{D}$ is a bounded plane domain with Lipschitz boundary. Moreover, $\breve{\Ga}_j$ is a connected set for every $j=0,\cdots, N$. In other words, $\{\breve{\Ga}_j: j=0,\cdots, N\}$ coincides with the family of all connected components of the set $P_+\cap \p\mc{D}$.
\end{description}

Hence ${\bf w}, {\bf k}$ and $p$ satisfy the following system in the plan domain $\mc{D}$:
\be\lab{a42}\begin{cases}
w_r\p_r w_r + w_z\p_z w_r -\f{w_{\th}^2}{r} + \p_r p =-\f{k_{\th}^2}{r},\\
w_r\p_r w_{\th}+ w_z\p_z w_{\th} +\f{w_r w_{\th}}{r} =0,\\
w_r\p_r w_z + w_z\p_z w_z + \p_z p=0,\\
w_r\p_r k_{\th}+ w_z\p_z k_{\th}-\f{w_r k_{\th}}{r}=0,\\
\p_r(rw_r) +\p_z (r w_z)=0.
\end{cases}\ee
These equations are satisfied for almost all $x\in \mc{D}$ and
\be\lab{a43}
{\bf w}(x)={\bf k}(x)=0\q \text{for $\mc{H}^1$-almost all $x\in P_+\cap\p \mc{D}$}.
\ee

We have the following integral estimates: ${\bf w}, {\bf k}\in W_{loc}^{1,2}(\mc{D})$,
\be\lab{a44}
\int_{\mc{D}}(|\nabla {\bf w}(r,z)|^2+|\na {\bf k}(r,z)|^2) r dr dz<\oo
\ee
and, by the Sobolev embedding theorem for three-dimensional domains, ${\bf w}, {\bf k}\in L^6(\Om)$, i.e.,
\be\lab{a45}
\int_{\mc{D}}(|{\bf w}(r,z)|^6+|{\bf k}(r,z)|^6) r dr dz <\oo.
\ee

Also, the condition $\nabla p\in L^{3/2}(\Om)$ can be written as
\be\lab{a46}
\int_{\mc{D}}|\nabla p(r,z)|^{3/2} r dr dz <\oo.
\ee
Denote by $\Phi= p +\f{|{\bf w}|^2}{2}+\f{|{\bf k}|^2}{2}$ the total head pressure corresponding to the solution $({\bf w}, {\bf k}, p)$. Obviously,
\be\lab{a411}
\int_{\mc{D}} r |\nabla \Phi(r,z)|^{3/2} dr dz<\oo.
\ee
Hence
\be\lab{a412}
\Phi\in W^{1,3/2}(\mc{D}_{\e})\q \q \forall \e>0.
\ee
We also have the important {\it Bernoulli's law}: for almost all $x\in\mc{D}$
\be\lab{bernoullilaw}
(w_r\p_r+ w_z\p_z)\Phi=0.
\ee

\subsection{Some results on Inviscid MHD equations.}\lab{inviscid mhd}

Since ${\bf w}, {\bf k}$ satisfy (\ref{o413}), then ${\bf w}={\bf k}\equiv 0$ on $\p\Om$ and $\na p\in L^{3/2}(\Om)$, then one can follow \cite{amick84} and \cite{kp83} to prove the following statement.
\bl\lab{al42}
{\it If ({\bf IMHD-AX}) are satisfied, then
\be\lab{a47}
\forall j\in \{0,1,\cdots, N\}\q \exists p_j\in \mbR: \q p(x)\equiv p_j\q \text{for $\mc{H}^2$-almost all $x\in\Ga_j$}.
\ee
In particular, by axial symmetry,
\be\lab{a48}
p(x)\equiv p_j\q \text{for $\mc{H}^1$-almost all $x\in \breve{\Ga}_j$.}
\ee
}
\el

We need a weak version of Bernoulli's law for a Sobolev solution $({\bf w}, {\bf k}, p)$ to the inviscid MHD equations (\ref{a42}).

From the last equality in (\ref{a42}) and from (\ref{a44}) it follows that there exists a stream function $\psi\in W_{loc}^{2,2}(\mc{D})$ such that
\be\lab{a49}
\f{\p\psi}{\p r}= -r w_z,\q \f{\p \psi}{\p z} = r w_r.
\ee

Fix a point $x_*\in \mc{D}$. For $\e>0$, denote by $\mc{D}_{\e}$ the connected component of $\mc{D}\cap \{(r,z): r>\e\}$ containing $x_*$. Since
\be\lab{a410}
\psi\in W^{2,2}(\mc{D}_{\e})\q\q \forall \e>0,
\ee
by Sobolev embedding theorem, $\psi\in C(\ol{\mc{D}_{\e}})$. Hence $\psi$ is continuous at points of $\ol{\mc{D}}\setminus O_z= \ol{\mc{D}}\setminus \{(0,z): z\in \mbR\}$. By the definition of $\psi$ and ${\bf w}={\bf k}\equiv 0$ on $\p\Om$, we see that all the boundary components are level sets of $\psi$.
\bl\lab{al421}
{\it If ({\bf IMHD-AX}) are satisfied, then there exist constants $\xi_0,\cdots, \xi_N\in \mbR$ such that $\psi(x)\equiv \xi_j$ on each curve $\breve{\Ga}_j$, $j=0,\cdots, N$.
}
\el

\bpf
In virtue of (\ref{a43}) and (\ref{a49}), we have $\nabla \psi(x)=0$ for $\mc{H}^1$-almost all $x\in \p\mc{D}\setminus O_z$. Then the Morse-Sard property (see \cite{bkk13}) implies that
\be\no
\textit{for any connected set } C\subset \p\mc{D}\setminus O_z,\q \exists \al=\al(C)\in \mbR: \psi(x)\equiv \al\q \forall x\in C.
\ee
Hence since $\breve{\Ga}_j$ are connected, the lemma follows.

\epf

By the properties of Sobolev functions ${\bf w}, {\bf k}, \psi,\Phi$ (see \cite{eg92}), we get the following
\bl\lab{al43}
{\it If conditions ({\bf IMHD-AX}) hold, then there exists a set $A_{{\bf w}}\subset \ol{\mc{D}}$ such that
\begin{enumerate} [(i)]
  \item $\mc{H}^1(A_{{\bf w}})=0$.
  \item For all $x=(r,z)\in \mc{D}\setminus A_{{\bf w}}$,
  \be\no
  \displaystyle\lim_{\rho\to 0} \fint_{B_{\rho}(x)} |{\bf w}(y)- {\bf w}(x)|^2 dy=\lim_{\rho\to 0} \fint_{B_{\rho}(x)} |{\bf k}(y)- {\bf k}(x)|^2 dy = \lim_{\rho\to 0} \fint_{B_{\rho}(x)} |\Phi(y)-\Phi(x)|^2 dy=0;
  \ee
  moreover, the function $\psi$ is differentiable at $x$, and
  \be\no
  \nabla \psi(x) = (-r w_z(x), r w_r(x)).
  \ee
  \item For every $\e>0$, there exists a set $U\subset \mbR^2$ with $\mc{H}_{\oo}^1(U)<\e, A_{{\bf w}}\subset U$, and such that the functions ${\bf w}, {\bf k}, \Phi$ are continuous on $\ol{\mc{D}}\setminus (U\cup O_z)$.
\end{enumerate}
}
\el

Then one can mimic the proof in \cite{kpr15as} to establish the following weak version of {\it Bernoulli's law}.
\bl\lab{al44b}({\it Bernoulli's Law}).
{\it Let conditions ({\bf IMHD-AX}) be valid, and let $A_{{\bf w}}$ be a set from Lemma \ref{al43}. For any compact connected set $K\subset \ol{\mc{D}}\setminus O_z$, the following property holds: if
\be\lab{a413}
\psi|_K= \text{const},
\ee
then
\be\lab{a414}
\Phi(x_1) =\Phi(x_2)\q \q \text{for all }x_1,x_2\in K\setminus A_{{\bf w}}.
\ee
}
\el
In particular, we can denote by $\Phi(K)$ the uniform constant $c\in \mbR$ such that $\Phi(x)=c$ for all $x\in K\setminus A_{{\bf w}}$, for any compact set $K\subset \ol{\mc{D}}\setminus O_z$ with $\psi_{K}=\textit{const}$. Moreover, $\Phi$ has some continuity properties when $K$ approaches the singularity axis $O_z$.

\bl\lab{al44}
{\it Assume that conditions ({\bf IMHD-AX}) are satisfied. Let $K_i$ be a sequence of compact sets with the following properties: $K_i\subset \ol{\mc{D}}\setminus O_z, \psi|_{K_i}=\text{const}$, and $\displaystyle\lim_{i\to \oo} \inf_{(r,z)\in K_i} r=0$, $\displaystyle\liminf_{i\to\oo} \sup_{(r,z)\in K_i} r>0$. Then $\Phi(K_i)\to p_0$ as $i\to \oo$.
}
\el

\bl\lab{al43b}
{\it If conditions ({\bf IMHD-AX}) are satisfied, then $p_0=\cdots= p_{M'}$, where $p_j$ are the constants from Lemma \ref{al42}.
}
\el
Heuristically, one can imagine that the axis $Oz$ is an ``almost" stream line, by Lemma \ref{al44b}, all the boundary components that intersects with the symmetry axis should share the same total head pressure $\Phi$, which immediately implies Lemma \ref{al43b}. Since the proof of Lemmas \ref{al42}-\ref{al43b} are quite similar to the proofs in \cite{kpr15as}, we omit the details.

\subsection{Obtaining a contradiction.} \lab{contradiction}

We consider three possible cases.
\begin{enumerate}[(a)]
  \item The maximum of $\Phi$ is attained on the boundary component intersecting the symmetry axis:
  \be\lab{a415}
  p_0 =\max_{j=0,\cdots,N} p_j =\sup_{x\in \Om} \Phi(x).
  \ee
  \item The maximum of $\Phi$ is attained on a boundary component that does not intersect the symmetry axis:
  \be\lab{a416}
  p_0< p_N= \max_{j=0,\cdots, N} p_j = \sup_{x\in \Om} \Phi(x).
  \ee
  \item The maximum of $\Phi$ is not attained on $\p\Om$:
  \be\lab{a417}
  \max_{j=0,\cdots, N} p_j < \sup_{x\in \Om} \Phi(x).
  \ee
\end{enumerate}

\subsubsection{The case $\sup_{x\in \Om} \Phi(x)= p_0$.} Adding a constant to the pressure $p$, we can assume that
\be\lab{a418}
p_0= \sup_{x\in\Om}\Phi(x)=0.
\ee

Since the identity $p_0=p_1=\cdots= p_N$ is impossible, we have that $p_j<0$ for some $j\in \{M'+1, N\}$. Recall that by Lemma \ref{al43b}, $p_0=p_1=\cdots =p_{M'}=0$. From equation ($\ref{o413}_1$) we obtain
\be\lab{a419}\begin{array}{ll}
0&= x\cdot\nabla p(x) +x\cdot({\bf w}(x)\cdot\nabla) {\bf w}(x)-x\cdot({\bf k}(x)\cdot\nabla) {\bf k}(x) \\
&=\text{div }[x p(x)+ ({\bf w}(x)\cdot x) {\bf w}(x)-({\bf k}(x)\cdot x){\bf k}(x)]- p(x) \text{div }x- |{\bf w}(x)|^2+|{\bf k}(x)|^2\\
&=\text{div }[x p(x)+ ({\bf w}(x)\cdot x) {\bf w}(x)-({\bf k}(x)\cdot x){\bf k}(x)]- 3\Phi(x) +\f 12|{\bf w}(x)|^2+\f{5}{2}|{\bf k}(x)|^2.
\end{array}\ee

Integrating it over $\p\Om$ and using (\ref{a418}), we derive a contradiction as follows
\be\no
0&\geq& \int_{\Om} [3\Phi(x)-\f 12 |{\bf w}(x)|^2-\f{5}{2}|{\bf k}(x)|^2] dx= \int_{\p\Om} p(x) (x\cdot {\bf n}) ds= \sum_{j=0}^N p_j \int_{\Ga_j}(x\cdot{\bf n}) ds\\\no
&=& \sum_{j=1}^N p_j \int_{\Om_j} \text{div }x dx= -3 \sum_{j=1}^N p_j |\Om_j|>0.
\ee
Hence we exclude the first case.

\subsubsection{The case $p_0<p_N=\sup_{x\in \ol{\Om}}\Phi(x)$.} We may assume that the maximum value is zero:
\be\lab{a420}
p_0<p_N =\max_{j=0,\cdots,N} p_j =\sup_{x\in\Om}\Phi(x)=0.
\ee
Then $p_0=\cdots =p_{M'}<0$.

Change (if necessary) the numbering of the boundary components $\Ga_{M'+1},\cdots, \Ga_{N-1}$ so that
\be\lab{a422}
{p}_j <0, j=0,\cdots, M, M\geq M',\\\lab{a423}
{p}_{M+1}=\cdots ={p}_N=0.
\ee

To remove a neighborhood of the singularity line $O_z$ from our consideration, we take $r_0>0$ such that the open set $\mc{D}_{\e}=\{(r,z)\in \mc{D}: r>\e\}$ is connected for every $\e\leq r_0$ (i.e. $\mc{D}_{\e}$ is a domain), and
\be\lab{a424}\begin{array}{ll}
\breve{\Ga}_j \subset \ol{\mc{D}_{r_0}}\q \text{and}\q \inf_{(r,z)\in \breve{\Ga}_j} r\geq 2r_0,\q j=M'+1,\cdots, N,\\
\breve{\Ga}_j\cap \ol{\mc{D}_{\e}}\q \text{is a connected set and }\sup_{(r,z)\in\breve{\Ga}_j\cap\ol{\mc{D}_{\e}}} r\geq 2r_0, j=0,\cdots, M', \e\in (0,r_0].
\end{array}\ee

Let a set $C\subset \ol{\mc{D}_{\e}}$ separate $\breve{\Ga}_i$ and $\breve{\Ga}_j$ in $\mc{D}_{\e}$ for some different indexes $i,j\in\{0,\cdots, N\}$; i.e., $\breve{\Ga}_i\cap \ol{\mc{D}_{\e}}$ and $\breve{\Ga}_j\cap \ol{\mc{D}_{\e}}$ lie in different connected components of $\ol{\mc{D}_{\e}}\setminus C$. Obviously, for $\e\in (0, r_0]$, there exists a constant $\de(\e)>0$ (not depending on $i, j, C$) such that the uniform estimate $\sup_{(r,z)\in C} r\geq \de(\e)$ holds. Moreover, the function $\de(\e)$ is nondecreasing. In particular,
\be\lab{a425}
\de(\e)\geq \de(r_0),\q \e\in (0, r_0].
\ee

In the following, we will construct an appropriate integration domain by using the level sets of $\Phi$ and $\Phi_n$. We need some information concerning the behavior of the limit total head pressure $\Phi$ on stream lines. Following \cite{kpr15annals} and \cite{kronrod50}, we introduce some facts of topology. By {\it continuum} we mean a compact connected set. We understand connectedness in the sense of general topology. A subset of a topological space is called {\it an arc} if it is homeomorphic to the unit interval $[0,1]$. Let $Q=[0,1]\times [0,1]$ be a square in $\mbR^2$, and let $f$ be a continuous function on $Q$. Denote by $E_t$ a level set of the function $f$, i.e., $E_t=\{x\in Q: f(x)=t\}$. A connected component $K$ of the level set $E_t$ containing a point $x_0$ is a maximal connected subset of $E_t$ containing $x_0$. By $T_f$ denote a family of all connected components of level sets of $f$.

We apply Kronrod's results to the stream function $\psi|_{\ol{\mc{D}_{\e}}}, \e\in (0, r_0]$. Accordingly, $T_{\psi, \e}$ means the corresponding Kronrod tree for the restriction $\psi|_{\ol{\mc{D}_{\e}}}$. Define the total head pressure on the Kronrod tree $T_{\psi,\e}$ as follows. Let $K\in T_{\psi, \e}$ with $\text{diam }K>0$. Take any $x\in K\setminus A_{{\bf w}}$ and put $\Phi(K)= \Phi(x)$. By the Bernoulli's Law in Lemma \ref{al44b}, the value $\Phi(K)$ is independent of the choice $x\in K\setminus A_{{\bf w}}$. Then $\Phi$ has the following continuity properties on stream lines.

\bl\lab{al45}(See Lemma 3.5 in \cite{kpr15annals}).
{\it Let $A, B\in T_{\psi,\e}$, where $\e\in (0,r_0]$, $\text{diam } A>0$, and $\text{diam }B>0$. Consider the corresponding arc $[A, B]\subset T_{\psi,\e}$ joining $A$ to $B$. Then the restriction $\Phi|_{[A, B]}$ is a continuous function.
}
\el

Denote by $B_0^{\e},\cdots, B_N^{\e}$ the elements of $T_{\psi, \e}$ such that $B_j^{\e}\supset \breve{\Ga}_j \cap \ol{\mc{D}_{\e}}$, $j=0,\cdots, M'$, and $B_j^{\e}\supset \breve{\Ga}_j, j=M'+1,\cdots, N$. By construction, $\Phi(B_j^{\e})<0$ for $j=0,\cdots, M$, and $\Phi(B_j^{\e})=0$ for $j=M+1,\cdots, N$. For $r>0$, let $L_r$ be the horizontal straight line $L_r=\{(r,z): z\in \mbR\}$. Then similar to Lemma 4.6 in \cite{kpr15annals}, we can find $r_*\in (0, r_0]$ and $C_j\in [B_j^{r_*}, B_N^{r_*}]$, $j=0,\cdots, M$, such that $\Phi(C_j)<0$ and $C\cap L_{r_*}= \emptyset$ for all $C\in [C_j, B_N^{r_*}]$.


We restrict our argument on the domain $\mc{D}_{r_*}$ and put $T_{\psi}= T_{\psi, r_*}$ and $B_j= B_j^{r_*}$. Since $\p\mc{D}_{r_*}\subset B_0\cup\cdots \cup B_N\cup L_{r_*}$ and the set $\{B_0,\cdots, B_N\}\subset T_{\psi}$ is finite, we can change $C_j$ (if necessary) such that
\be\lab{a426}
&&\forall j=0,\cdots, M,\q C_j\in [B_j, B_N],\q \Phi(C_j)<0,\\\lab{a427}
&&C\cap \p\mc{D}_{r_*}=\emptyset\q \q \forall C\in [C_j, B_N).
\ee

Observe that $\Ga_j\cap L_{r_*}\neq \emptyset$ for $j=0,\cdots, M'$. Therefore, if a cycle $C\in T_{\psi}$ separates $\Ga_N$ from $\Ga_0$ and $C\cap \p\mc{D}_{r_*}=\emptyset$, then $C$ separates $\Ga_N$ from $\Ga_j$ for all $j=1,\cdots, M'$. So we can take $C_0= \cdots =C_{M'}$ and consider only the Kronrod arcs $[C_{M'}, B_N],\cdots, [C_M, B_N]$. Recall that a set $\mc{Z}\subset T_{\psi}$ has $T$-measure zero if $\mc{H}^1(\{\psi(C): C\in \mc{Z}\})=0$.

\bl\lab{al47}
{\it For every $j=M',\cdots, M$, $T$-almost all $C\in [C_j, B_N]$ are $C^1$-curves homeomorphic to the circle. Moreover, there exists a subsequence $\Phi_{n_l}$ such that the sequence $\Phi_{n_l}|_{C}$ converges to $\Phi|_C$ uniformly $\Phi_n|_C \rightrightarrows \Phi|_C$ on $T$-almost all cycles $C\in [C_j, B_N]$.
}\el
%
%

Without loss of generality, we assume that the subsequence $\Phi_{n_l}$ coincides with $\Phi_n$. Besides, cycles satisfying the assertion of Lemma \ref{al47} will be called {\it regular cycles}. From Lemmas \ref{al47} and Lemma 3.6 in \cite{kpr15annals}, we can conclude that
\be\label{lemma}
\mc{H}^1(\{\Phi(C): C\in [C_j, B_N]\q \textit{and $C$ is not a regular cycle}\})=0, j=M',\cdots, M.
\ee

Setting $\al=\displaystyle\max_{j=M',\cdots, M}\min_{C\in [C_j, B_N]} \Phi(C)$, then by (\ref{a426}), $\al<0$. By (\ref{lemma}), we can find a sequence of positive values $t_i\in (0,-\al), i\in \mb{N}$ with $t_{i+1}= \f 12 t_i$, such that the implication
\be\no
\Phi(C)=-t_i\Rightarrow \textit{$C$ is a regular cycle}
\ee
holds for every $j=M',\cdots, M$ and for all $C\in [C_j, B_N]$. Consider the natural order on the arc $[C_j, B_N]$, namely, $C'<C^{''}$ if $C^{''}$ is closer to $B_N$ than $C'$. For $j=M',\cdots, M$ and $i\in\mb{N}$, put
\be\no
A_i^j=\max\{C\in [C_j, B_N]: \Phi(C)=-t_i\}.
\ee
Then each $A_i^j$ is a regular cycle and $A_i^j\subset \mc{D}_{r_*}$. In particular, for each $i\in \mb{N}$, the compact set $\displaystyle\cup_{j=M'}^M A_i^j$ is separated from $\p \mc{D}_{r_*}$ and $\textit{dist}(\displaystyle\cup_{j=M'}^M A_i^j,\p \mc{D}_{r_*})>0$. Then for each $i$ and for sufficiently small $h>0$, we have the inclusion $\{x\in \mc{D}_{r_*}: \textit{dist}(x, \breve{\Ga}_N)<h\}\subset \mc{D}_{r_*}\setminus (\displaystyle\cup_{j=M'}^{M} A_i^j)$. Denote by $V_i$ the connected component of the open set $\mc{D}_{r_*}\setminus (\displaystyle\cup_{j=M'}^{M} A_i^j)$ which encloses the set $\{x\in \mc{D}_{r_*}: \textit{dist}(x, \breve{\Ga}_N)<h\}$. Then we have
\be\no
\{x\in \mc{D}_{r_*}: \textit{dist}(x, \breve{\Ga}_N)<h\}\cap \p V_i= A_i^{M'}\cup\cdots \cup A_i^{M}.
\ee
By the construction, the sequence of domains $V_i$ is decreasing; i.e. $V_i\supset V_{i+1}$. Hence the sequence of sets $(\p \mc{D}_{r_*})\cap (\p V_i)$ is nonincreasing. Every set $(\p \mc{D}_{r_*})\cap (\p V_i)$ consists of several components $\breve{\Ga}_l$ with $l>M$. Since there are only finitely many components $\Ga_l$, then we can conclude that for sufficiently large $i$, the set $(\p \mc{D}_{r_*})\cap (\p V_i)$ is independent of $i$. So we can assume that $(\p \mc{D}_{r_*})\cap (\p V_i)=\breve{\Ga}_K\cup\cdots \cup \breve{\Ga}_N$, where $K\in \{M+1,\cdots, N\}$. Hence,
\be\lab{a428}
\p V_i= A_i^{M'}\cup\cdots \cup A_i^M\cup \breve{\Ga}_K\cup\cdots \cup \breve{\Ga}_N.
\ee

By Lemma \ref{al47}, we have the uniform convergence $\Phi_n|_{A_i^j}\rightrightarrows \Phi(A_i^j)$ as $n\to \oo$, then for each $i\in \mb{N}$ there exists $n_i$ such that for all $n\geq n_i$
\be\no
\Phi_n|_{A_i^j}<-\f{7}{8} t_i,\q \Phi_n|_{A_{i+1}^j}>-\f 58t_i\q \forall j =M',\cdots, M.
\ee
Then
\be\no
\forall t\in[\f 58 t_i,\f 78t_i]\q \forall n\geq n_i\q \Phi_n|_{A_i^j}<-t,\q \Phi_n|_{A_{i+1}^j}>-t\q \forall j=M',\cdots, M.
\ee

Accordingly, for $n\geq n_i$ and $t\in [\f58 t_i,\f 78t_i]$, we can define $W_{in}^j(t)$ as the connected component of the open set
$\{x\in V_i\setminus \ol{V_{i+1}}: \Phi_n(x)>-t\}$ with $\p W_{in}^j(t)\supset A_{i+1}^j$ and put
\be\no
W_{in}(t)= \bigcup_{j=M'}^M W_{in}^j(t),\q S_{in}(t)= (\p W_{in}(t))\cap (V_i\setminus \ol{V_{i+1}}).
\ee

By construction, $\Phi_n\equiv -t$ on $S_{in}(t)$ and
\be\no
\p W_{in}(t)= S_{in}(t)\cup A_{i+1}^{M'}\cup\cdots \cup A_{i+1}^M,
\ee
and the set $S_{in}(t)$ separates $A_i^{M'}\cup\cdots \cup A_i^M$ from $A_{i+1}^{M'}\cup\cdots \cup A_{i+1}^M$. Since $\Phi_n\in W^{2,2}_{loc}(\Om)$, by the Morse-Sard theorem, for almost all $t\in [\f 58 t_i,\f78 t_i]$, the level set $S_{in}(t)$ consists of finitely many $C^1$-cycles and $\Phi_n$ is differentiable in classical sense at every point $x\in S_{in}(t)$ with $\na \Phi_n(x)\neq 0$. We will say the values $t\in [\f58 t_i,\f78 t_i]$ having the above property are $(n,i)$-regular. Therefore, $\wt{S_{in}}(t)$ is a finite union of smooth surfaces (tori), and by construction,
\be\lab{a429}
\int_{\wt{S_{in}}(t)}\nabla \Phi_n \cdot {\bf n} d S=- \int_{\wt{S_{in}}(t)} |\nabla \Phi_n| dS<0,
\ee
where ${\bf n}$ is the unit outward normal vector to $\p \wt{W_{in}}(t)$.

For $h>0$, denote $\Ga_h=\{x\in\Om: \text{dist}(x, \Ga_K\cup\cdots\cup \Ga_N)=h\}$, $\Om_h=\{x\in\Om: \text{dist}(x,\Ga_K\cup\cdots \cup \Ga_N)<h\}$. Since the distance function $\text{dist}(x,\p\Om)$ is $C^1$-regular and the norm of its gradient is equal to one in the neighborhood of $\p\Om$, there is a constant $\de_0>0$ such that for every $h\leq \de_0$, the set $\Ga_h$ is a union of $N-K+1$ $C^1$-smooth surfaces homeomorphic to the torus, and
\be\lab{a430}
\mc{H}^2(\Ga_h)\leq c_0\q \forall h\in (0,\de_0],
\ee
where the constant $c_0= 3\mc{H}^2(\Ga_K\cup\cdots\cup \Ga_N)$ is independent of $h$.

%
%
%

\bl\lab{al49}
{\it For any $i\in \mb{N}$, there exists $n(i)\in \mb{N}$ such that for every $n\geq n(i)$ and for almost all $t\in [\f 58 t_i,\f78t_i]$, the inequality
\be\lab{a432}
\int_{\wt{S_{in}}(t)} |\na \Phi_n| dS\leq \mc{F} t
\ee
holds with the constant $\mc{F}$ independent of $t, n$ and $i$.
}\el

\bpf
By a direct calculation, (\ref{a42}) implies
\be\lab{a431}\begin{array}{ll}
\na \Phi&= \na\f{1}{2}|{\bf w}|^2- ({\bf w}\cdot\na){\bf w}+ \na\f{1}{2}|{\bf k}|^2 + ({\bf k}\cdot\na){\bf k}\\
&=[\na {\bf w}-(\na {\bf w})^T]\cdot {\bf w} + [\na {\bf k}+(\na {\bf k})^T]\cdot {\bf k}.
\end{array}\ee
Since $\Phi\neq \textit{const}$ on $\wt{V_i}$, (\ref{a431}) implies $\int_{\wt{V_i}} |\na {\bf w}-(\na {\bf w})^T|^2+ |\na {\bf k}+(\na {\bf k})^T|^2 dx>0$ for every $i$. Hence, from the weak convergence $\na{\bf u}_n \rightharpoonup \na {\bf w}$ and $\na{\bf h}_n \rightharpoonup \na {\bf k}$ in $L^2(\Om)$ it follows that for any $i\in \mb{N}$, there exist constants $\e_i>0, \de_i\in (0,\de_0)$ and $k_i'\in\mb{N}$ such that
\be\no
\ol{\Om_{\de_i}}\cap \wt{A_i^j} = \ol{\Om_{\de_i}} \cap \wt{A_{i+1}^j} =\emptyset, \q j=M',\cdots, M,
\ee
and for all $n\geq n_i'$
\be\lab{a431'}
\int_{\wt{V_{i+1}}\setminus \Om_{\de_i}} \left(|\na {\bf u}_n-(\na {\bf u}_n)^T|^2+ |\na {\bf h}_n+(\na {\bf h}_n)^T|^2\right) dx >\e_i.
\ee

Fix $i\in \mb{N}$. We assume that $n\geq n_i$. Since we have removed a neighborhood of the singularity line $O_z$, we can use the Sobolev embedding theorem in the {\it plane} domain $\mc{D}_{r_*}$. The uniformly boundedness of $\|\Phi_n\|_{W^{1,3/2}(\mc{D}_{r_*})}$ imply that the norm $ \|\Phi_n\|_{L^6(\mc{D}_{r_*})}$ and then $\|\Phi_n\na \Phi_n\|_{L^{6/5}(\mc{D}_{r_*})}$ are also uniformly bounded. Finally we have
\be\lab{a433}
\|\Phi_n\na \Phi_n\|_{L^{6/5}(\wt{\mc{D}_{r_*}})}\leq \textit{const}.
\ee

Fix a sufficiently small $\si>0$ (the exact value of $\si$ will be specified below), and take the parameter $\de_{\si}\in (0, \de_i]$ small enough to satisfy the following conditions:
\be\lab{a434}
&&\Om_{\de_{\si}} \cap \wt{A_i^j} =\Om_{\de_{\si}} \cap \wt{A_{i+1}^j} =\emptyset,\q j=M',\cdots, M,\\\lab{a435}
&&\int_{\Ga_h} \Phi_n^2 dS <\si^2\q \forall h\in (0,\de_{\si}] \q \forall n\geq n'.
\ee
The last estimate follows from the identity $\Phi|_{\Ga_K\cup\cdots\cup\Ga_N}\equiv 0$, the weak convergence $\Phi_n \rightharpoonup \Phi$ in the space $W^{1,3/2}(\Om)$, and (\ref{a433}).

By a direct calculation, (\ref{a310}) implies
\be\no
\nabla \Phi_n &=& -\nu_n \text{curl }\text{curl }{\bf u}_n + [\na {\bf u}_n- (\na {\bf u}_n)^T]\cdot {\bf u}_n \\\no
&\q&+[\na {\bf h}_n+ (\na {\bf h}_n)^T]\cdot {\bf h}_n+\na\times {\bf f}_n.
\ee

Then using the Stokes theorem, we obtain
\be\no
\int_{S} \nabla \Phi_n\cdot {\bf n} d S&=& \int_{S}\left([\na {\bf u}_n- (\na {\bf u}_n)^T]\cdot {\bf u}_n\right)\cdot {\bf n} dS+\int_{S}\left([\na {\bf h}_n+ (\na {\bf h}_n)^T]\cdot {\bf h}_n\right)\cdot {\bf n} dS.
\ee

Now, fix a sufficiently small $\e>0$. The exact value of $\e$ will be specified below. For a given sufficiently large $n\geq n'$, we follow the argument in Lemma 3.8 of \cite{kpr15annals} to find a number $\ol{h_n}\in (0,\de_{\si})$ such that the estimates
\be\lab{a436}
\b|\int_{\Ga_{\ol{h_n}}} \na\Phi_n\cdot {\bf n} dS\b| &\leq& 2 \int_{\Ga_{\ol{h_n}}}(|{\bf u}_n|\cdot |\na {\bf u}_n|+|{\bf h}_n|\cdot |\na {\bf h}_n|) d S<\e,\\\lab{a437}
\int_{\Ga_{\ol{h_n}}} (|{\bf u}_n|^2+|{\bf h}_n|^2) d S&\leq& C_{\e} \nu_n^2
\ee
hold, where $C_{\e}$ is independent of $n$ and $\si$.

Now, for $(n,i)$-regular value $t\in [\f 58 t_i,\f 78 t_i]$, consider the domain
\be\no
\Om_{i\ol{h_{n}}}(t)= \wt{W_{in}(t)}\cup (\ol{\wt{V_{i+1}}}\setminus \ol{\Om_{\ol{h_n}}}).
\ee

By construction, $\p \Om_{i\ol{h_n}}(t)= \Ga_{\ol{h_n}}\cup \wt{S_{in}}(t)$. Also using (\ref{a310}), we know
\be\no
\De\Phi_n&=&\Delta p_n + |\na {\bf u}_n|^2 + |\na {\bf h}_n|^2+ {\bf u}_n\cdot \Delta {\bf u}+ {\bf h}_n\cdot \De {\bf h}_n \\\no
&=&-\text{div}(( {\bf u}_n\cdot \na) {\bf u}_n)+ \text{div}(( {\bf h}_n\cdot\na) {\bf h}_n)+ |\na  {\bf u}_n|^2+ |\na  {\bf h}_n|^2-\f1{\nu_n}\b((\na\times {\bf f}_n)\cdot {\bf u}_n+ (\na\times {\bf g}_n)\cdot {\bf h}_n\b)\\\no
&\q&+\f1{\nu_n}\b(({\bf u}_n\cdot\na) \f{|{\bf u}_n|^2}{2}+ {\bf u}_n\cdot\na p_n- {\bf u}_n\cdot(({\bf h}_n\cdot\na) {\bf h}_n)\b)+\f1{\nu_n}\b(({\bf u}_n\cdot\na) \f{|{\bf h}_n|^2}{2}- {\bf h}_n\cdot(({\bf h}_n\cdot\na) {\bf u}_n)\b)\\\no
&=&- \sum_{i,j=1}^3 \p_i u_{nj}\p_j u_{ni} +|\na {\bf u}_n|^2+ |\na {\bf h}_n|^2+ \sum_{i,j=1}^3 \p_{i} h_{nj} \p_j h_{ni}+\f{1}{\nu_n}({\bf u}_n\cdot\na) \Phi_n\\\no
&\q&-\f{1}{\nu_n}({\bf h}_n\cdot\na) ({\bf u}_n\cdot {\bf h}_n)- \f{1}{\nu_n} \b((\na\times {\bf f}_n)\cdot {\bf u}_n+ (\na\times {\bf g}_n)\cdot {\bf h}_n\b)\\\no
&=&\f 1{\nu_n}\text{div }(\Phi_n {\bf u}_n)+\f{1}{2}|\na {\bf u}_n-(\na {\bf u}_n)^T|^2+\f{1}{2}|\na {\bf h}_n+ (\na {\bf h}_n)^T|^2\\\lab{a348}
&\q&- \f{1}{\nu_n} \b((\na\times {\bf f}_n)\cdot {\bf u}_n+ (\na\times {\bf g}_n)\cdot {\bf h}_n\b),
\ee
where we have used the special structure of ${\bf u}_n$ and ${\bf h}_n$, so that $({\bf h}_n\cdot\na)({\bf u}_n\cdot {\bf h}_n)\equiv 0$. Integrating the equation (\ref{a348}) over the domain $\Om_{i\ol{h_n}}(t)$, we obtain
\be\no
&\quad& \int_{\wt{S_{in}}}\na \Phi_n\cdot {\bf n} ds + \int_{\Ga_{\ol{h_n}}}\na \Phi_n\cdot {\bf n} ds\\\no
&=& \int_{\Om_{i\ol{h_n}}(t)} \f{1}{2}|\na {\bf u}_n-(\na {\bf u}_n)^T|^2+\f{1}{2}|\na {\bf h}_n+ (\na {\bf h}_n)^T|^2 dx -\f 1{\nu_n} \int_{\Om_{i\ol{h_n}}(t)}\left((\na\times {\bf f}_n)\cdot {\bf u}_n+ (\na\times {\bf g}_n)\cdot {\bf h}_n\right)dx\\\no
&\q&+ \f1{\nu_n} \int_{\wt{S_{in}}}\Phi_n {\bf u}_n\cdot {\bf n} ds +\f1{\nu_n}\int_{\Ga_{\ol{h_n}}}\Phi_n {\bf u}_n\cdot {\bf n} ds \\\no
&=& \int_{\Om_{i\ol{h_n}}(t)} \f{1}{2}|\na {\bf u}_n-(\na {\bf u}_n)^T|^2+\f{1}{2}|\na {\bf h}_n+ (\na {\bf h}_n)^T|^2 dx -\f 1{\nu_n} \int_{\Om_{i\ol{h_n}}(t)}\left((\na\times {\bf f}_n)\cdot {\bf u}_n+ (\na\times {\bf g}_n)\cdot {\bf h}_n\right)dx\\\no
&\q&+ \f1{\nu_n}\int_{\Ga_{\ol{h_n}}}\Phi_n {\bf u}_n\cdot {\bf n} ds- t \la_n \ol{\mc{F}},
\ee
where $\ol{\mc{F}}= (\mc{F}_{M'}+\cdots + \mc{F}_M)$. In view of (\ref{a436}), we can estimate
\be\no
\int_{\wt{S_{in}}} |\na \Phi_n| ds &\leq & t\mc{F} +\e+\f 1{\nu_n} \int_{\Om_{i\ol{h_n}}(t)}\left((\na\times {\bf f}_n)\cdot {\bf u}_n+ (\na\times {\bf g}_n)\cdot {\bf h}_n\right)dx  \\\no
&\q&-\int_{\Om_{i\ol{h_n}}(t)} \b(\f{1}{2}\b|\na {\bf u}_n-(\na {\bf u}_n)^T\b|^2+\f{1}{2}\b|\na {\bf h}_n+ (\na {\bf h}_n)^T\b|^2\b) dx\\\lab{a350}
&\quad& + \f{1}{\nu_n} \b(\int_{\Ga_{\ol{h_n}}} \Phi_n^2 ds\b)^{\f 12} \b(\int_{\Ga_{\ol{h_n}}} |{\bf u}_n|^2 ds\b)^{\f 12},
\ee
with $\mc{F}= |\ol{\mc{F}}|$. By definition, $\f{1}{\nu_n}\|\na\times {\bf f}_n\|_{L^2(\Om)}= \la_n\nu_n\|\na\times {\bf f}\|_{L^2(\Om)}\to 0$ as $n\to \oo$. Therefore,
\be\no
\b|\f 1{\nu_n} \int_{\Om_{i\ol{h_n}}(t)}\left((\na\times{\bf f}_n)\cdot {\bf u}_n+(\na\times{\bf g}_n)\cdot {\bf h}_n\right) dx\b|\leq \e
\ee
for sufficiently large $n$. Using inequalities (\ref{a435}) and (\ref{a437}) in (\ref{a350}), we obtain
\be\no
\int_{\wt{S_{in}}} |\na \Phi_n| ds &\leq &t \mc{F}+ 2\e + \si \sqrt{C_{\e}}- \int_{\Om_{i\ol{h_n}}(t)} \f{1}{2}|\na {\bf u}_n-(\na {\bf u}_n)^T|^2+\f{1}{2}|\na {\bf h}_n+ (\na {\bf h}_n)^T|^2  dx \\\no
&\leq& t\mc{F} + 2\e + \si\sqrt{C_{\e}}- \int_{\wt{V_{i+1}}\setminus \Om_{\de_i}} \f{1}{2}|\na {\bf u}_n-(\na {\bf u}_n)^T|^2+\f{1}{2}|\na {\bf h}_n+ (\na {\bf h}_n)^T|^2  dx,
\ee
where $C_{\e}$ is independent of $n$ and $\si$. Choosing $\e=\f 16 \e_i$, $\si=\f{\e_i}{3\sqrt{C_{\e}}}$, and a sufficiently large $n$, from (\ref{a431'}) we obtain $2\e+ \si\sqrt{C_{\e}}- \int_{\wt{V_{i+1}}\setminus \Om_{\de_i}} \f{1}{2}|\na {\bf u}_n-(\na {\bf u}_n)^T|^2+\f{1}{2}|\na {\bf h}_n+ (\na {\bf h}_n)^T|^2 dx\leq 0$. We have finished the proof.

\epf

Now we can derive a contradiction by using the Co-area formula.
\bl\lab{al39}
{\it Assume that $\Om\subset \mbR^3$ is a bounded domain of type (\ref{a21}) with $C^2$-smooth boundary $\p\Om$, $(\na\times {\bf f},\na\times {\bf g})\in W_{AS}^{1,2}(\Om)\times W_{ASoS}^{1,2}(\Om)$, and $({\bf a}, {\bf b})\in W_{AS}^{3/2,2}(\p\Om)\times W_{ASoS}^{3/2,2}(\p\Om)$ satisfies (\ref{comp1})-(\ref{comp2}). Then assumptions ({\bf MHD-AX}) and (\ref{a420}) lead to a contradiction.
}\el

The proof of Lemma \ref{al39} can be obtained by slightly modifying the proof of Lemma 3.9 of \cite{kpr15annals}, i.e., replacing Hausdorff measure $\mc{H}^1$ by $\mc{H}^2$, and the curves $S_{in}(t)$ by the surfaces $\wt{S_{in}}(t)$ in the corresponding integrals and the details are omitted. Therefore, we have excluded the second case.


\subsubsection{The case $\sup_{x\in\Om} \Phi(x)> \max_{j=0,\cdots,N} p_j$.}

Assume that (\ref{a417}) is satisfied, and set $\si= \max_{j=0,\cdots,N} p_j$. Then we can find a compact connected set $F\subset \mc{D}\setminus A_{{\bf w}}$ such that $\text{diam}(F)>0$, $\psi|_F=\textit{const}$, and $\Phi(F)>\si$. We may assume that $\si<0$ and $\Phi(F)=0$. We still need to separate $F$ from $\p\mc{D}$ by regular cycles, and take a number $r_0>0$ such that $F\subset \mc{D}_{r_0}$, the open set $\mc{D}_{\e}= \{(r,z)\in \mc{D}: r>\e\}$ is connected for every $\e\leq r_0$, and conditions (\ref{a424}) are satisfied. Then for $\e\in (0, r_0]$, we can consider the behavior of $\Phi$ on the Kronrod trees $T_{\psi,\e}$ corresponding to the restrictions $\psi|_{\ol{\mc{D}_{\e}}}$. Denote by $F^{\e}$ the element of $T_{\psi,\e}$ containing $F$. Using the same procedure as previous, we can find $r_*\in (0, r_0]$ and $C_j\in [B_j^{r_*}, F^{r_*}]$, $j=0,\cdots, N$, such that $\Phi(C_j)<0$ and $C\cap L_{r_*}=\emptyset$ for all $C\in [C_j, F^{r_*}]$.

Set $T_{\psi}= T_{\psi, r_*}, F^*= F^{r_*}$, and $B_j= B_j^{r_*}$, i.e. $B_j\in T_{\psi}$ and $B_j\supset \breve{\Ga}_j \cap \ol{\mc{D}_{r_*}}$. As above, we can change $C_j$ so that
\be\no
&&\q\forall j=0,\cdots, N\q C_j\in [B_j, F^*],\q \Phi(C_j)<0,\\\no
&&C\cap \p\mc{D}_{r_*}=\emptyset\q \forall C\in [C_j, F^*],\q \textit{and}\q C_0=\cdots =C_{M'}.
\ee


Similarly, we should construct an appropriate integration domain by using the level sets of $\Phi$ and $\Phi_n$. Take positive numbers $t_i= 2^{-i} t_0$, regular cycles $A_i^j\in [C_j, F^*]$ with $\Phi(A_i^j)= -t_i$, and the set $S_{in}(t)$ with $\Phi_n|_{S_{in}(t)}\equiv -t$ separating $A_i^{M'}\cup\cdots \cup A_i^N$ from $A_{i+1}^{M'}\cup\cdots \cup A_{i+1}^N$, etc. Argued as in Lemma \ref{al49} and \ref{al39}, we can derive a similar contradiction as before. Therefore, we have finished the proof of Theorem \ref{at41}.

{\bf Acknowledgement.}  Weng's research was supported by Basic Science Research Program through the National Research Foundation of Korea(NRF) funded by the Ministry of Education, Science and Technology (2015049582). The author would like to thank Prof. Dongho Chae and Prof. Zhouping Xin for their interests in this work and constant encouragement and support. The author want to thank the referees for their careful reading and important suggestion and improvements.

\end{document}